\documentclass[a4paper,12pt,final]{amsart}

\usepackage{amssymb}
\usepackage[OT4]{fontenc}
\usepackage{texdraw}

\theoremstyle{definition}
\newtheorem{definition}{Definition}
\theoremstyle{plain}
\newtheorem{theorem}[definition]{Theorem}
\newtheorem{lemma}[definition]{Lemma}
\newtheorem{conjecture}[definition]{Conjecture}

\newtheorem{proposition}[definition]{Proposition}
\theoremstyle{remark}
\newtheorem{remark}[definition]{Remark}
\newtheorem{example}[definition]{Example}

\DeclareMathOperator{\edim}{edim}
\DeclareMathOperator{\vdim}{vdim}
\DeclareMathOperator{\rank}{rank}

\DeclareMathOperator{\charfield}{char}
\DeclareMathOperator{\bdeg}{bdeg}
\DeclareMathOperator{\sgn}{sgn}

\DeclareMathOperator{\EoLS}{EoLS}

\def\field{\mathbb{K}}
\def\N{\mathbb{N}}
\def\Z{\mathbb{Z}}
\def\R{\mathbb{R}}
\def\proj{\mathbb{P}}
\def\sys{\mathcal{L}}
\def\tto{\longrightarrow}
\def\up{\uparrow}
\def\Mat{\mathcal{M}}

\begin{document}

\title[Reduction method for linear systems of plane curves...]{Reduction method for linear systems of plane curves with
base fat points}

\author{Marcin Dumnicki}








\begin{abstract}
In the paper we develop a new method of proving non-speciality
of a linear system with base fat points in general position.
Using this method we
show that the Hirschowitz-Harbourne Conjecture holds
for systems with base points of equal multiplicity bounded by $42$.
\end{abstract}

\maketitle

\section{Introduction}

Let $\field$ be a field of characteristics zero, 
$\N=\{0,1,2,\dots\}$, $\N^\ast=\{1,2,\dots\}$, let $\R$
denote the field of real numbers.

\begin{definition}
Let $D\subset\mathbb N^2$ be finite (any such set will be called a
\textit{diagram}), let $m_{1},\dots,m_{r} \in \N^{*}$, let $p_{1}, \dots, p_{r}
\in \field^{2}$. Define the $\field$-vector space $\sys_{D}(m_{1}p_{1},\dots,
m_{r}p_{r}) \subset \field[X,Y]$:
\begin{align*}
\sys_{D}(m_{1}p_{1}&,\dots,m_{r}p_{r}) := \\
\big\{ & f = \sum_{(\beta_{1},\beta_{2}) \in D} 
c_{(\beta,\beta_{2})}X^{\beta_{1}}Y^{\beta_{2}} \mid 
c_{(\beta_{1},\beta_{2})} \in \field,
\frac{\partial^{\alpha_{1}+\alpha_{2}} f}{\partial X^{\alpha_{1}}\partial Y^{\alpha_{2}}}(p_{j})
= 0, \\ & \alpha_{1}+\alpha_{2} < m_{j}, j=1,\dots,r \big\}.
\end{align*}
Define the \textit{dimension} of the \textit{system of curves
$\sys_{D}(m_{1},\dots,m_{r})$} to be
$$\dim \sys_{D}(m_{1},\dots,m_{r}) := \min_{\{p_{j}\}_{j=1}^{r}, p_{j} \in \field^{2}} \dim_{\field}
\sys_{D}(m_{1}p_{1}, \dots,m_{r}p_{r}) - 1.$$
\end{definition}

\begin{remark}
If points $p_{1},\dots,p_{r}$ are in general position we have
$$\dim \sys_{D}(m_{1},\dots,m_{r}) = \dim_{\field} \sys_{D}(m_{1}p_{1},\dots,m_{r}p_{r}) - 1.$$
\end{remark}

The system $\sys_{D}(m_{1},\dots,m_{r})$ can be understood as a vector space
of curves generated by monomials with exponents from $D$ having multiplicities at least
$m_{1},\dots,m_{r}$ in $r$ general points.

\begin{definition}
Let $L=\sys_{D}(m_{1},\dots,m_{r})$ be a system of curves. Define
the \textit{virtual dimension $\vdim$ of $L$} and the \textit{expected dimension $\edim$ of $L$}
\begin{align*}
\vdim L & := \# D - 1 - \sum_{j=1}^{r} \binom{m_{j}+1}{2}, \\
\edim L & := \max \{ \vdim L, -1\}.
\end{align*}
\end{definition}

Intuitively, we should have $\dim L = \edim L$.

\begin{definition}
We say that a \textit{system of curves $L$ is special}
if
$$\dim L > \edim L.$$
Otherwise we say that \textit{$L$ is non-special}.
\end{definition}

Observe that by linear algebra we have always $\dim L \geq \edim L$
since multiplicity $m$ imposes $\binom{m+1}{2}$ conditions.

\section{The Hirschowitz--Harbourne Conjecture}

For systems of the form $\sys_{d}(m_{1},\dots,m_{r}) :=
\sys_{D}(m_{1},\dots,m_{r})$, $D=\{\alpha \in \N^{2} \mid |\alpha| \leq d\}$
the well-known
Hirschowitz--Harbourne Conjecture giving geometrical
description of the speciality of a system  was formulated in \cite{HCON}.
To formulate this conjecture consider the blowing-up
$\pi : \widetilde\proj^{2} \to \proj^{2}$ in $r$ general points with
exceptional divisors $E_{1},\dots,E_{r}$.

\begin{definition}
A curve $C \subset \proj^{2}$ is said to be a \textit{$-1$-curve}
if it is irreducible and the self-intersection $\widetilde{C}^{2}$ of its proper
transform $\widetilde{C} \subset \widetilde\proj^{2}$ is equal to $-1$.
\end{definition}

\begin{conjecture}[Hirschowitz--Harbourne]\label{hirsch}
A system $L = \sys_{d}(m_{1},\dots,m_{r})$ is special if and only if
there exists a $-1$-curve $C \subset \proj^{2}$ such that
$$\widetilde{L}.\widetilde{C} \leq -2,$$
where $\widetilde{L} := |d\pi^{*}(\mathcal O_{\mathbb P^2}(1)) - \sum_{j=1}^{r}m_{j}E_{j}|$.
\end{conjecture}

This conjecture was studied by many authors, we refer only to the recent
results. For homogenous systems ($m_1=m_{2}=m_{3}=\dots=m_{r}=:m$), 
the above conjecture holds for $m \leq 20$ (see \cite{CCMO,CMir}). 
In the general case the Conjecture holds for multiplicities bounded by $11$
(see \cite{hiha11}). Both these results were obtained with the help of computers.

For further information about the above conjecture see for example
\cite{CIL}.

In this paper we develop a new method (called ``diagram cutting'') based
on some properties of matrices assigned to systems of curves.
This method allows providing easy proofs
of non-speciality for large family of systems. Moreover, these proofs
can often be found argorithmically with a computer program.
Sometimes this needs a lot of computations, but then it can easily be checked
``by hand''.

As a result of the method we show that in order to check non-speciality
of all homogeneous systems of bounded multiplicity $m$ it is enough to
check a finite number of cases. This result was obtained in purely theoretical way.

The second result is Theorem \ref{hiha42} stating that 
Conjecture \ref{hirsch} holds for homogeneous multiplicities bounded by $42$.
This result was obtained using a computer program.

\section{Diagram cutting method\label{secmet}}

Before introducing the method we must establish the notation and
say when a system is non-special in the language of matrices.

\begin{definition}
Let $j \in \N^{*}$, $\alpha \in \N^{2}$. Define the mapping
$$
\varphi_{j,\alpha}:\field[X,Y] \ni f\mapsto 
\frac{\partial^{|\alpha|} f}{\partial X^{\alpha}}(P_{j,X},P_{j,Y})
\in\field[P_{j,X},P_{j,Y}],
$$
where $P_{j,X},P_{j,Y}$ are new indeterminates used instead of $X,Y$.
\end{definition}

We will use the following notation: $\Mat(n,k;R)$ will denote
the set of all matrices with $n$ rows, $k$ columns, and coefficients
from $R$ (a ring or a field). For an $M \in \Mat(n,k;R)$ we write
$M_{[j,\ell]}$ for the element of $M$ in the $j$-th row and the $\ell$-th column.

\begin{definition}
Let $L=\sys_{D}(m_{1},\dots,m_{r})$ be a system of curves, let
$D = \{(\alpha_{1,X},\alpha_{1,Y}),\dots,(\alpha_{n,X},\alpha_{n,Y})\}$,
$\alpha_{i,X},\alpha_{i,Y} \in \N$ for $i=1,\dots,n$.
Let $\mathcal A=\{(j,\beta)\in \N \times \N^{2} \mid|\beta|<m_j,\;j=1,\dots,r\} 
= \{\mathfrak{a}_{1},\dots,\mathfrak{a}_{c}\}$.
Define the matrix
$M(L) \in \Mat(c,n;\field[P_{i,X},P_{i,Y}]_{i=1}^{r})$ as
$$M(L)_{[j,k]} := \varphi_{\mathfrak{a}_{j}}(X^{\alpha_{k,X}}Y^{\alpha_{k,Y}}).$$
\end{definition}

For given points $p_{1}=(p_{1,X},p_{1,Y}),\dots,p_{r}=(p_{r,X},p_{r,Y}) \in
\field^{2}$ we will use the natural evaluation mapping
$$\nu_{p_{1},\dots,p_{r}} : \field[P_{i,X},P_{i,Y}]_{i=1}^{r} 
\ni f \mapsto f|_{P_{i,X} \mapsto p_{i,X}, P_{i,Y} \mapsto p_{i,Y}}
\in \field.$$

\begin{proposition}
Let $L=\sys_{D}(m_{1},\dots,m_{r})$ be a system of curves.
Then $\dim L = \# D - \rank M(L) - 1$.
\end{proposition}

\begin{proof}
Let $p_{1},\dots,p_{r} \in \field^{2}$ be points in general position.
Consider the linear mapping
$$
\Phi : \big\{ f=\sum_{(\alpha_{X},\alpha_{Y}) \in D}
c_{(\alpha_{X},\alpha_{Y})}X^{\alpha_{X}}Y^{\alpha_{Y}} \big\} \ni f 
\mapsto (\nu_{p_{1},\dots,p_{r}} \circ \varphi_{\mathfrak{a}_{j}}(f))_{j=1}^{c} 
\in \field^{c}.
$$
We have $\sys_{D}(m_{1}p_{1},\dots,m_{r}p_{r}) = \ker \Phi$.
Let $M$ denote the matrix of $\Phi$ in bases
$\{(\alpha_{1,X},\alpha_{1,Y}),\dots,(\alpha_{n,X},\alpha_{n,Y})\}$,
$\{\mathfrak{a}_{1},\dots,\mathfrak{a}_{c}\}$. We have
$$\nu_{p_{1},\dots,p_{r}}(M(L)_{[j,k]}) = M_{[j,k]},$$
hence $\rank M(L) = \rank M$ (we use
the facts that $\charfield \field = 0$ and $p_{1},\dots,p_{r}$ are in
general position). Now we compute
$$\dim L = n - \rank M - 1 = \# D - \rank M(L) - 1.$$
\end{proof}

\begin{definition}
Define the bidegree $\bdeg : \field[P_{i,X},P_{i,Y}]_{i=1}^{r} \to \N^{2}$
given by
$\bdeg(P_{i,X}) := (1,0), \bdeg(P_{i,Y}) := (0,1)$ for $i=1,\dots,r$.
\end{definition}

\begin{proposition}
Let $\sys_{D}(m_{1},\dots,m_{r})$ be a system of curves.
Let $M$ be a square submatrix of $M(L)$ of size $s \in \N^{*}$.
After a necessary renumbering of columns and rows we can assume
that $M$ is given by columns
$(\alpha_{1,X},\alpha_{1,Y}),\dots,(\alpha_{s,X},\alpha_{s,Y})$ and
rows $\mathfrak{a}_{1},\dots,\mathfrak{a}_{s}$. 
Then $\det M$ is a bihomogeneous (w.r.t. $\bdeg$) polynomial
of bidegree $(\sum_{i=1}^{s} \alpha_{i,X},\sum_{i=1}^{s} \alpha_{i,Y})
- \gamma$, where $\gamma \in \N^{2}$ depends
only on the choice of rows.
\end{proposition}

\begin{proof}
$$\det M = \sum_{\sigma \in S_{s}} \sgn(\sigma) M_{[1,\sigma(1)]} \dots M_{[s,\sigma(s)]}.$$
For $M_{[j,k]} \neq 0$ we have 
$\bdeg M_{[j,k]} = (\alpha_{k,X},\alpha_{k,Y}) - \beta_{j}$,
where $\mathfrak{a}_{j} = (\ell_{j},\beta_{j})$ for some $\ell_{j} \in \N$, $\beta_{j} \in \N^{2}$.
Hence
$$
\bdeg M_{[1,\sigma(1)]} \cdots M_{[s,\sigma(s)]} = 
\bigg(\sum_{i=1}^{s} \alpha_{i,X},\sum_{i=1}^{s} \alpha_{i,Y}\bigg)
- \sum_{i=1}^{s} \beta_{i}.
$$
We finish the proof by taking $\gamma = \sum_{i=1}^{s}\beta_{i}$.
\end{proof}

\begin{proposition}
\label{onemult}
Let $m \in \N^{*}$, let $D \subset \N^{2}$. Assume that
$\# D = \binom{m+1}{2}$. Then
$L=\sys_{D}(m)$ is non-special if and only if $D$ do not lie on a curve of
degree $m-1$.
\end{proposition}

\begin{proof}
From the previous proof we can see that $\det M(L) = c \cdot f$, 
where $f$ is a monomial, $c \in \field$.
Let
$D = \{(\alpha_{1,X},\alpha_{1,Y}),\dots,(\alpha_{n,X},\alpha_{n,Y})\}$.
For $\beta=(\beta_{X},\beta_{Y})$, $|\beta| < m$ we have
\begin{multline*}
M(L)_{[(1,\beta), j]} = \varphi_{(1,\beta)}(X^{\alpha_{j,X}}Y^{\alpha_{j,Y}}) \\ 
= \prod_{k=1}^{\beta_{X}} (\alpha_{j,X}-k+1) \cdot
\prod_{k=1}^{\beta_{Y}} (\alpha_{j,Y}-k+1) \cdot
P_{1,X}^{\alpha_{j,X} - \beta_{X}}P_{1,Y}^{\alpha_{j,Y} - \beta_{Y}}.
\end{multline*}
Since we are only interested in the value of $c$ we would like to compute
the determinant of $M = M(L)_{X \mapsto 1, Y \mapsto 1}$.
By a suitable sequence of row operations we can
change $M$ into the following matrix $M'$
$$M'_{[(1,\beta),j]} = \alpha_{j,X}^{\beta_{X}}\alpha_{j,Y}^{\beta_{Y}},
\qquad \det M \neq 0 \iff \det M' \neq 0.$$
Now $\det M' = 0$ if and only if rows of $M'$ are linearly dependent, but
this happens if and only if $D$ lies on a curve of degree $m-1$.
\end{proof}

Now we can present the diagram cutting method and prove that it can be used
to bound the dimension of a system of curves.

\begin{theorem}
\label{linecut}
Let $m_{1},\dots,m_{r},m_{r+1},\dots,m_{p} \in \N^{*}$, let $D \subset \N^{2}$ be a diagram,
let $F : \N^{2} \ni (\alpha_{1},\alpha_{2}) \mapsto r_{1}\alpha_{1}+r_{2}\alpha_{2}+r_{0} \in \R$
be an affine function, $r_{0},r_{1},r_{2} \in \R$. Let
\begin{align*}
D_{1} & := \{ (\alpha_{1},\alpha_{2}) \in D \mid F(\alpha_{1},\alpha_{2}) < 0 \}, \\
D_{2} & := \{ (\alpha_{1},\alpha_{2}) \in D \mid F(\alpha_{1},\alpha_{2}) > 0 \}.
\end{align*}
If $L_{2}:=\sys_{D_{2}}(m_{r+1},\dots,m_{p})$ is non-special and $\vdim L_{2} = -1$
then 
$$\dim \sys_{D}(m_{1},\dots,m_{p}) \leq \dim \sys_{D_{1}}(m_{1},\dots,m_{r}).$$
\end{theorem}

\begin{proof}
Put $L_{1}:=\sys_{D_{1}}(m_{1},\dots,m_{r})$.
We can compute the dimension of a system $L:=\sys_{D}(m_{1},\dots,m_{p})$ as
$\dim L = \# D - \rank M(L) - 1$.
As $D=D_{1} \cup D_{2}$, in an appropriate basis the matrix $M(L)$ is of the following form
$$M(L) = \left[ \begin{array}{c|c}
M(L_{1}) & K_{1} \\ \hline
K_{2} & M(L_{2}) \\ \end{array} \right].$$
Take a maximal non-zero minor $M'$ of $M(L_{1})$ and consider the following
square submatrix of $M(L)$:
$$M = \left[ \begin{array}{c|c}
M' & K_{1}' \\ \hline
K_{2}' & M(L_{2}) \\ \end{array} \right],$$
where $K_{1}'$ and $K_{2}'$ are suitable submatrices of $K_{1}$ and $K_{2}$.
It suffices to show that $\det M \neq 0$.
The columns of $M'$ are indexed by elements from some $D_{1}' \subset D_{1}$
hence the columns of $M$ are indexed by $D':=D_{1}' \cup D_{2}$.
Let $U = [M' \; K_{1}']$, $L = [K_{2}' \; M(L_{2})]$ be submatrices
of $M$, let 
$$\mathcal{C} := \{ C \subset D' \mid \#C = \#D_{2} \}.$$
For $C \in \mathcal{C}$ define $L_{C}$ (respectively $U_{C}$)
as the submatrix of $L$ (resp. $U$) given by taking the columns indexed by
elements of $C$.
Now we can compute
$$\det M = \sum_{C \in \mathcal{C}} \epsilon(C) \det L_{C} \det U_{D' \setminus C},$$
for $\epsilon(C) = \pm 1$.
Observe that $\det L_{C} \in \field[P_{r+1,X},P_{r+1,Y},\dots,P_{p,X},P_{p,Y}]$,
$\det U_{D' \setminus C} \in \field[P_{1,X},P_{1,Y},\dots,P_{r,X},P_{r,Y}]$ and
consider $\det M$ as a polynomial of
indeterminates $\{P_{\ell,X},P_{\ell,Y}\}_{\ell > r}$ over
$\field[\{P_{\ell,X},P_{\ell,Y}\}_{\ell \leq r}]$.
We have
$$\det M = \pm \det M' \det M(L_{2}) + f.$$
Assume that $\det M=0$. As $\det M' \neq 0$ and $\det M(L_{2}) \neq 0$ from
the assumptions, we must have 
another non-zero term 
$g \in \field[\{P_{\ell,X},P_{\ell,Y}\}_{\ell > r}]$ appearing
in $f$ such that $\bdeg(g) = \bdeg(\det M(L_{2}))$. The term
$g$ was given by some $C \in \mathcal{C}$, $C \neq D_{2}$.
Since 
$$
\bdeg(g) = \sum_{(\alpha_{1},\alpha_{2}) \in C} (\alpha_{1},\alpha_{2}) - \gamma, 
\quad
\bdeg(\det M(L_{2}))  = \sum_{(\alpha_{1},\alpha_{2}) \in D_{2}} (\alpha_{1},\alpha_{2}) - \gamma
$$
we have
$$F \bigg( \sum_{(\alpha_{1},\alpha_{2}) \in C} (\alpha_{1},\alpha_{2}) \bigg) = 
F \bigg( \sum_{(\alpha_{1},\alpha_{2}) \in D_{2}} (\alpha_{1},\alpha_{2}) \bigg).$$
As $F$ is an affine form and $\# C = \# D_{2}$ we have
$$\sum_{(\alpha_{1},\alpha_{2}) \in C} F(\alpha_{1},\alpha_{2}) = 
\sum_{(\alpha_{1},\alpha_{2}) \in D_{2}} F(\alpha_{1},\alpha_{2}),$$
but from the definition of $D_{2}$ it is possible only when $C = D_{2}$,
contradiction.
\end{proof}

\section{Reduction of homogeneous systems}

We will use the following notation for a sequence of multiplicities:

\begin{definition}
Let $m_{1},\dots,m_{r} \in \N^{*}$, $k_{1},\dots,k_{r} \in \N$. Define
$$(m_{1}^{\times k_{1}},\dots,m_{r}^{\times k_{r}}) =
(\underbrace{m_{1},\dots,m_{1}}_{k_{1}},\dots,\underbrace{m_{r},\dots,m_{r}}_{k_{r}}).$$
\end{definition}

We will use diagrams of the following form:

\begin{definition}
Let $a_{1},\dots,a_{n},u_{1},\dots,u_{n} \in \N$. Define
$$(a_{1}^{\up u_{1}},\dots,a_{n}^{\up u_{n}}) :=
\bigcup_{i=1}^{n} \{i-1\} \times \{u_{i},\dots,u_{i}+a_{i}-1\} \subset \N^{2}.$$
\end{definition}

\begin{example}
$$
\begin{array}{rcl}
\begin{texdraw}
\drawdim pt
\linewd 0.5
\move(0 0)
\fcir f:0 r:0.2
\move(0 8)
\fcir f:0 r:0.2
\move(0 16)
\fcir f:0 r:0.2
\move(0 24)
\fcir f:0 r:0.2
\move(0 32)
\fcir f:0 r:0.2
\move(8 0)
\fcir f:0 r:0.2
\move(8 8)
\fcir f:0 r:0.2
\move(8 16)
\fcir f:0 r:0.2
\move(8 24)
\fcir f:0 r:0.2
\move(8 32)
\fcir f:0 r:0.2
\arrowheadtype t:V
\move(0 0)
\avec(32 0)
\move(0 0)
\avec(0 56)
\htext(37 0){$\mathbb{N}$}
\htext(-13 40){$\mathbb{N}$}
\move(0 24)
\fcir f:0 r:1.2
\move(0 32)
\fcir f:0 r:1.2
\move(8 0)
\fcir f:0 r:1.2
\textref h:C v:C
\htext(16 -15){Fig. 1. Diagram $(2^{\up 3},1^{\up 0})$}
\end{texdraw}
&
\hskip3cm & 
\begin{texdraw}
\drawdim pt
\linewd 0.5
\move(0 0)
\fcir f:0 r:0.2
\move(0 8)
\fcir f:0 r:0.2
\move(0 16)
\fcir f:0 r:0.2
\move(0 24)
\fcir f:0 r:0.2
\move(0 32)
\fcir f:0 r:0.2
\move(8 0)
\fcir f:0 r:0.2
\move(8 8)
\fcir f:0 r:0.2
\move(8 16)
\fcir f:0 r:0.2
\move(8 24)
\fcir f:0 r:0.2
\move(8 32)
\fcir f:0 r:0.2
\move(16 0)
\fcir f:0 r:0.2
\move(16 8)
\fcir f:0 r:0.2
\move(16 16)
\fcir f:0 r:0.2
\move(16 24)
\fcir f:0 r:0.2
\move(16 32)
\fcir f:0 r:0.2
\move(24 0)
\fcir f:0 r:0.2
\move(24 8)
\fcir f:0 r:0.2
\move(24 16)
\fcir f:0 r:0.2
\move(24 24)
\fcir f:0 r:0.2
\move(24 32)
\fcir f:0 r:0.2
\arrowheadtype t:V
\move(0 0)
\avec(48 0)
\move(0 0)
\avec(0 56)
\htext(53 0){$\mathbb{N}$}
\htext(-13 40){$\mathbb{N}$}
\move(0 24)
\fcir f:0 r:1.2
\move(0 32)
\fcir f:0 r:1.2
\move(8 0)
\fcir f:0 r:1.2
\move(24 16)
\fcir f:0 r:1.2
\move(24 24)
\fcir f:0 r:1.2
\move(24 32)
\fcir f:0 r:1.2
\textref h:C v:C
\htext(24 -15){Fig. 2. Diagram $(2^{\up 3},1^{\up 0}, 0^{\up 0}, 3^{\up 2})$}
\end{texdraw}
 \\
\end{array}
$$
\end{example}

Observe that $\# (a_{1}^{\up u_{1}},\dots,a_{n}^{\up u_{n}}) = \sum_{i=1}^{n} a_{i}$.

\begin{definition}
Let $\Z$ denote the set of all integers.
We say that two diagrams $D_{1}$, $D_{2}$ are
\textit{equivalent} (or $D_{1}$ \textit{is equivalent to} $D_{2}$) 
if there exists $\alpha \in \Z^{2}$ such that
$D_{1} = D_{2} + \alpha$.
\end{definition}

\begin{remark}
Observe that diagram $(0^{\up 0},\dots,0^{\up 0},a_{1}^{\up u_{1}},\dots,a_{n}^{\up u_{n}})$
is equivalent to $(a_{1}^{\up u_{1}},\dots,a_{n}^{\up u_{n}})$ is equivalent
to $(a_{1}^{\up u_{1}+u},\dots,a_{n}^{\up u_{n}+u})$.
\end{remark}

\begin{proposition}
Let $m_{1},\dots,m_{r} \in \N^{*}$, let $D_{1}$, $D_{2}$ be diagrams.
If $D_{1}$ and $D_{2}$ are equivalent, then
$\dim \sys_{D_{1}}(m_{1},\dots,m_{r}) = \dim \sys_{D_{2}}(m_{1},\dots,m_{r})$.
\end{proposition}

\begin{proof}
Let $D_{1} + (\alpha_{1},\alpha_{2}) = D_{2}$, let $p_{1},\dots,p_{r} \in \field^{2}$ be
points in general position. The maps
\begin{align*}
\sys_{D_{1}}(m_{1}p_{1},\dots,m_{r}p_{r}) \ni f & \mapsto X^{\alpha_{1}}Y^{\alpha_{2}}f \in \sys_{D_{2}}(m_{1}p_{1},\dots,m_{r}p_{r}), \\
\sys_{D_{2}}(m_{1}p_{1},\dots,m_{r}p_{r}) \ni f & \mapsto X^{-\alpha_{1}}Y^{-\alpha_{2}}f \in \sys_{D_{1}}(m_{1}p_{1},\dots,m_{r}p_{r}) \\
\end{align*}
are well-defined (we can assume that none of the coordinates of
$p_{1},\dots,p_{r}$ are zero), linear and inverse to each other.
\end{proof}

\begin{lemma}
\label{backtriangle}
Let $m \in \N^{*}$, $D=(1^{\up m-1},2^{\up m-2},\dots,(m-1)^{\up 1},m^{\up 0})$.
Then $\sys_{D}(m)$ is non-special and $\vdim \sys_{D}(m) = -1$.
\end{lemma}

\begin{proof}
By Lemma \ref{onemult} it is enough to show that $D$ 
(Fig. 3 shows an example for $m=3$) do not lie on a curve $C$
of degree $m-1$. Let $L_{j} = \{ (x,y) \in \R^{2} \mid x+y+j=0 \}$, $j=0,\dots,m-1$.
Observe that $\# (D \cap L_{j}) = j+1$ so by Bezout Theorem and induction
we have $\bigcup_{j=0}^{m-1}L_{j} \subset C$, a contradiction.
\end{proof}

\begin{remark}
Observe that we can do the same for diagram 
$(m^{\up 0},(m-1)^{\up 0},\dots,1^{\up 0})$.
\end{remark}

\begin{lemma}
\label{twotriangles}
Let $m \in \N^{*}$, $D=(m^{\up m},m^{\up m-1},\dots,m^{\up 0})$.
Then $\sys_{D}(m^{\times 2})$ is non-special and $\vdim \sys_{D}(m^{\times 2}) = -1$.
\end{lemma}

\begin{proof}
Let $F = y-m+\frac{1}{2}$. Observe that $D=D_{1} \cup D_{2}$
(Fig. 4 shows an example for $m=3$), where
\begin{align*}
D_{1} & := \{ p \in D \mid F(p) < 0 \} = (0,1^{\up m-1},2^{\up m-2},\dots,m-1^{\up 1},m^{\up 0}), \\
D_{2} & := \{ p \in D \mid F(p) > 0 \} = (m^{\up m},m-1^{\up m},\dots,1^{\up m}).\\
\end{align*}
The diagram $D_{1}$ is equivalent to 
$(1^{\up m-1},2^{\up m-2},\dots,(m-1)^{\up 1},m^{\up 0})$ hence from Lemma
\ref{backtriangle} the system $\sys_{D_{1}}(m)$ is non-special.
The diagram $D_{2}$ is equivalent to
$(m^{\up 0},(m-1)^{\up 0},\dots,1^{\up 0})$ so $\sys_{D_{2}}(m)$ is
non-special. As $\# D_{2} = \binom{m+1}{2}$, we can use
Theorem \ref{linecut} to obtain non-speciality of $\sys_{D}(m^{\times 2})$.
\end{proof}

\begin{lemma}
\label{singlelayer}
Let $m,k \in \N^{*}$, $D=(m^{\up k-1},m^{\up k-2},\dots,m^{\up 0})$.
If $m+1|k$ then $L=\sys_{D}(m^{\times 2\frac{k}{m+1}})$ is non-special
and $\vdim L = -1$.
\end{lemma}

\begin{proof}
We proceed by induction on $k$. For $k=m+1$ we use the previous Lemma.
Let $k > m+1$. Take $F=x-(m+1)+\frac{1}{2}$. Observe that $D=D_{1} \cup D_{2}$ 
(Fig. 5 shows an example for $m=3$, $k=12$), where
\begin{align*}
D_{1} & := \{ p \in D \mid F(p) < 0 \} = (m^{\up k-1},\dots,m^{\up k-1-m}), \\
D_{2} & := \{ p \in D \mid F(p) > 0 \} = (0,\dots,0,m^{\up k-(m+1)-1},\dots,m^{\up 0}). \\
\end{align*}
The diagram $D_{1}$ is equivalent to $(m^{\up m},\dots,m^{\up 0})$, hence
from Lemma \ref{twotriangles} the system $\sys_{D_{1}}(m^{\times 2})$
is non-special.
The diagram $D_{2}$ is equivalent to
$(m^{\up k-(m+1)-1},\dots,m^{\up 0})$ and from the induction hypothesis
we get that system $\sys_{D_{2}}(m^{\times 2\frac{k-(m+1)}{m+1}})$ is non-special.
Now, Theorem \ref{linecut} allows us to finish the proof.
\end{proof}

\begin{lemma}
\label{fatlayer}
Let $m, k, h \in \N^{*}$, $D=(h^{\up k-1},h^{\up k-2},\dots,h^{\up 0})$.
If $m+1|k$ and $m|h$ then $L=\sys_{D}(m^{\times 2\frac{kh}{m(m+1)}})$ is non-special
and $\vdim L = -1$.
\end{lemma}

\begin{proof}
We proceed by induction on $h$. If $h=m$ then the situation was treated
in the previous Lemma. Let $h > m$.
Take $F=y+x-(k-1+m)+\frac{1}{2}$. Observe that $D=D_{1} \cup D_{2}$
(Fig. 6 shows an example for $m=3$, $k=12$, $h=9$), where
\begin{align*}
D_{1} & := \{ p \in D \mid F(p) < 0 \} = (m^{\up k-1},m^{\up k-2},\dots,m^{\up 0}), \\
D_{2} & := \{ p \in D \mid F(p) > 0 \} = ((h-m)^{\up k-1+m},\dots,(h-m)^{\up m}). \\
\end{align*}
According to Lemma \ref{singlelayer} the system
$\sys_{D_{1}}(m^{\times 2\frac{k}{m+1}})$ is non-special.
The diagram $D_{2}$ is equivalent to
$((h-m)^{\up k-1},(h-m)^{\up k-2},\dots,(h-m)^{\up 0})$ and by the induction
hypothesis the system
$\sys_{D_{2}}(m^{\times 2\frac{k(h-m)}{m(m+1}})$ is non-special.
Again we finish the proof by using Theorem \ref{linecut}.
\end{proof}

\begin{definition}
Let $m \in \N^{*}$, $h=m(m+1)$.
Define the set (called the \textit{end of layer systems})
\begin{align*}
\EoLS(m) = \{ & \sys_{D}(m^{\times 2k+h-1}) \mid \\
& D = (h^{\up k-1},\dots,h^{\up 0},(h-1)^{\up 0},\dots,1^{\up 0}),
k=1,\dots,m+1 \}.
\end{align*}
Observe that for every $L \in \EoLS(m)$ we have $\vdim L = -1$.
\end{definition}

\begin{lemma}
\label{fulllayer}
Let $m,k \in \N^{*}$, $h=m(m+1)$, $p=2k+h-1$, 
$D=(h^{\up k-1},h^{\up k-2},\dots,h^{\up 0},(h-1)^{\up 0},\dots,1^{\up 0})$.
If the set $\EoLS(m)$ contains only non-special systems then
the system $L=\sys_{D}(m^{\times p})$ is non-special.
\end{lemma}

\begin{proof}
Take $k_{1}, k_{2} \in \N$ such that
$k=k_{1}(m+1)+k_{2}$, $1 \leq k_{2} \leq m+1$. Put $F=x-k_{1}(m+1)+\frac{1}{2}$.
Observe that $D=D_{1} \cup D_{2}$ (Fig. 7 shows an example for $m=3$, $k=11$), where
\begin{align*}
D_{1} & := \{ p \in D \mid F(p) < 0 \} = (h^{\up k-1},h^{\up k-2},\dots,h^{\up k_{2}}), \\
D_{2} & := \{ p \in D \mid F(p) > 0 \} = (h^{\up k_{2}-1},\dots,h^{\up 0},h-1^{\up 0},\dots,1^{\up 0}). \\
\end{align*}
The diagram $D_{1}$ is equivalent to the diagram
$(h^{\up k-k_{2}-1},h^{\up k-k_{2}-2},\dots,h^{\up 0})$ and since
$m+1|k-k_{2}$, it follows from Lemma \ref{fatlayer} that
$\sys_{D_{1}}(m^{\times 2\frac{(k-k_{2})h}{m(m+1)}})$ is non-special.
The system $\sys_{D_{2}}(m^{\times 2k_{2}+h-1}) \in \EoLS(m)$ and
we can use Theorem \ref{linecut} to complete the proof.
\end{proof}

\begin{theorem}
\label{finitely}
Let $m, d_{L} \in \N^{*}$. Assume that for $d=d_{L},\dots,d_{L}+m(m+1)$ every
system $\sys_{d}(m^{\times p})$, $p \in \N$, is non-special. Moreover,
assume that the set $\EoLS(m)$ contains only non-special systems. Then
for any $d \geq d_{L}, p \in \N$ the system
$\sys_{d}(m^{\times p})$ is non-special.
\end{theorem}

\begin{proof}
We proceed by induction on $d$.
For $d_{L} \leq d \leq d_{L} + m(m+1)$ the proof is obvious. Take 
$d > d_{L}+m(m+1)$. 
We want to show that the system $\sys_{D}(m^{\times p})$ is
non-special, where $D=((d+1)^{\up 0},\dots,1^{\up 0})$. 
Take $h=m(m+1)$, $F=y+x-(d-h)+\frac{1}{2}$. Observe that $D=D_{1} \cup D_{2}$ 
(Fig. 8 shows an example for $m=3$, $d_{L}=3$, $d=16$), where
\begin{align*}
D_{1} & := \{ p \in D \mid F(p) < 0 \} = (d+1-h^{\up 0},\dots,1^{\up 0}), \\
D_{2} & := \{ p \in D \mid F(p) > 0 \} = (h^{\up d+1-h},\dots,h^{0},h-1^{\up 0},\dots,1^{\up 0}). \\
\end{align*}
As $d+1-h \geq d_{L}$ we may use the induction hypothesis for
the system $\sys_{D_{1}}(m^{\times p-(2d-h+3)})$.
By Lemma \ref{fulllayer} the system $\sys_{D_{2}}(m^{\times 2d-h+3})$ is
non-special. Again we finish the proof by using Theorem \ref{linecut}.
\end{proof}

\begin{definition}
Let $m,d_{0} \in \N^{*}$. Put
\begin{multline*}
S(m,d_{0}) := \EoLS(m) \cup \\
\{ \sys_{d}(m^{\times r}) \mid
\vdim \sys_{d}(m^{\times r}) \geq -2m^{2}, d_{0} \leq d \leq d_{0}+m(m+1),
r \in \N^{*} \}.
\end{multline*}
\end{definition}

\begin{theorem}
\label{setsm}
Let $m,d_{0} \in \N$. If the set $S(m,d_{0})$ contains only non-special systems
then every system $\sys_{d}(m^{\times r})$ for $d \geq d_{0}$, $r \in \N$
is non-special.
\end{theorem}

\begin{proof}
By Theorem \ref{finitely} it suffices to 
show that every system $\sys_{d}(m^{\times r})$, $r \in \N$,
$d=d_{0},\dots,d_{0}+m(m+1)$, is non-special.
Let $L=\sys_{d}(m^{\times r})$. If $\vdim L \geq -2m^{2}$ then $L$ is non-special
by the assumptions. Let $\vdim L < -2m^{2}$. From now on
we use the notations and theory of reductions
introduced in \cite{Dum,hiha11}. We want to apply a sequence of $r$ weak $m$-reductions to
the diagram $D = \{ \alpha \in \N^{2} \mid |\alpha| \leq d \}$ to end with the
empty diagram,
$$D \stackrel{m\textrm{ w}}{\tto} D_{1} \stackrel{m\textrm{ w}}{\tto}
D_{2} \stackrel{m\textrm{ w}}{\tto} D_{3} \stackrel{m\textrm{ w}}{\tto} \dots
\stackrel{m\textrm{ w}}{\tto} \varnothing.$$
Following the notations of \cite{hiha11} let us consider a diagram
$D=(a_{1},a_{2},\dots,a_{n})$. Observe that the $m$-reduction is not possible
only if $a_{i}=a_{i+1}<m$. As $D$ is the result of a sequence of weak
$m$-reductions this can only happen for $i \leq 2m$. While performing a $m$-weak
reduction we use at most $m$ additional points for
each $a_{i}$, $i=1,\dots,2m$, and for each $i$ it is sufficient to do it only once.
So we use at most $2m^{2}$ additional points to reduce $D$ to the empty
diagram, hence if $\vdim L < -2m^{2}$ then $L$ is non-special.
\end{proof}

\begin{example}
$$

$$
\end{example}

\section{Homogeneous systems with bounded multiplicity}

\begin{theorem}
\label{hiha42}
The Hirschowitz-Harbourne Conjecture holds for homogeneous systems with multiplicities
bounded by $42$.
\end{theorem}

\begin{proof}
For $m < 20$ the result can be found in \cite{CMir}.
For $m=20,\dots,42$ we choose $d_{0}(m)=3m$.
In order to
check the conjecture we have to do the following:
\begin{enumerate}
\item Find all non-special systems among $\sys_{d}(m^{\times r})$ for 
$d \leq d_{0}(m)$ (there are only finitely many of them). Next, for every such system
we must show that it satisfies the Hirschowitz-Harbourne conjecture. This was done
with the help of computer programs. By the proof of Theorem \ref{setsm}
the maximal size of matrix (for $m=42$) can be $8128 \times 11656$, but
in most cases the combination of reduction method and Cremona transformation
gives an immediate answer.
\item
For every system from the set $S(m,d_{0}(m))$ we must prove its
non-speciality.
\end{enumerate}
As the set $S(m,d_{0}(m))$ contains systems with big size of diagrams, this
cannot be done without preparations.
For a system $L=\sys_{D}(m^{\times r}) \in \EoLS(m)$ we use reduction method described in
\cite{hiha11} to reduce the problem to the question about non-speciality
of $L'=\sys_{D'}(m^{\times 5})$ for some diagram $D'$. 
For $m=42$ this forces us to compute the determinant
of $4515 \times 4515$ matrix $43$ times --- this can be done.
For the rest of the systems from $S(m,d_{0}(m))$ we use the following fact
(see Proposition 28 in \cite{hiha11}).

\begin{theorem}
Let $m_{1},\dots,m_{r} \in \N^{*}$. There exists a diagram $D$ with
the property: if $\sys_{D}(m_{1},\dots,m_{r})$ is non-special then for all $d \in \N$
the system $\sys_{d}(m_{1},\dots,m_{r})$ is non-special.
\end{theorem}

So for each $r$ such that $\sys_{d}(m^{\times r}) \in S(m, d_{0}(m))$ we have to check
only one system $\sys_{D}(m^{\times r})$ for some diagram $D$ depending on $r$.
We also reduce this system to $\sys_{D'}(m^{\times 9})$.
In \cite{MYWWW} the reader can find the table with the actual number of cases to be checked,
as well as all necessary software with instructions, on how to perform the tests.
\end{proof}

\begin{remark}
The test decribed above can also be performed for greater values of $m$,
but for each $m \geq 43$ this will take at least several days of computation.
It seems that now it would be better to reorganize the method.
\end{remark}

\section{Closing remarks}

\begin{remark}
There exists another method of proving non-speciality of given system
(or for a family of systems) based on blowing-up the projective space
introduced by C. Ciliberto and R. Miranda (\cite{CMir}). It seems
that diagram cutting method is different and sometimes works better.
Moreover, all definitions and results of section \ref{secmet} can
be easily carried over to the higher-dimensional case of the
systems of polynomials in $n$ variables vanishing (with multiplicities)
at points in general position. This is not known for the method
of C. Ciliberto and R. Miranda.
\end{remark}

\begin{remark}
Observe that Theorem \ref{finitely} can be reformulated to the following.
\begin{theorem}
If the set $\EoLS(m)$ contains only non-special systems and
$\sys_{d}(m_{1},\dots,m_{r})$ is non-special then
$$\sys_{d+m(m+1)}(m_{1},\dots,m_{r},m^{\times p})$$
is non-special, where $p=2d+m(m+1)+1$.
\end{theorem}
This shows that in order to find all non-special systems of the form
$\sys_{d}(m_{1},\dots,m_{r})$ with $m_{i} \leq M$, $i=1,\dots,r$ it
is sufficient to check a finite number of cases.
\end{remark}

\begin{example}
We will show non-speciality of the system $L=\sys_{21}(7^{\times 6}, 6^{\times 4},1)$
by diagram cutting method. The proof (found by a computer) can be easily
read off from the picture. The system $L$ was studied in \cite{Yang}.
$$

$$
\end{example}

\medskip

\noindent
\textsc{
Marcin Dumnicki \\
Institute of Mathematics, Jagiellonian University, \\
Reymonta 4, 30-059 Krak\'ow, Poland \\
}

\end{document}